\documentclass[oneside,english]{amsart}
\usepackage[T1]{fontenc}
\usepackage[latin1]{inputenc}
\usepackage{color}
\usepackage{amsthm}
\usepackage{amssymb}
\usepackage[all]{xy}

\makeatletter
\numberwithin{equation}{section} 
\numberwithin{figure}{section} 
\theoremstyle{plain}

\@ifundefined{definecolor}
 {\@ifundefined{definecolor}
 {\usepackage{color}}{}
}{}

\usepackage{amsthm}

\usepackage[all]{xy}

\numberwithin{equation}{section} 
\numberwithin{figure}{section} 
\theoremstyle{plain}

\theoremstyle{plain}
 \newtheorem{thm}{Theorem}[section]
  \newtheorem{prop}[thm]{Proposition}
   \newtheorem{fact}[thm]{Fact}
   \newtheorem{cor}[thm]{Corollary}
  \newtheorem{rem}[thm]{Remark}
  \newtheorem{lem}[thm]{Lemma}
   \newtheorem*{claim}{Claim}
   \theoremstyle{definition}
   \newtheorem{defn}[thm]{Definition}
  \theoremstyle{remark}
  \newtheorem{notation}[thm]{Notation}
\usepackage[all]{xy}

\newcommand{\nc}{\newcommand}
\nc{\Q}{\ensuremath{\mathbb{Q}}\xspace}
\renewcommand{\phi}{\varphi}
\nc{\acl}{\operatorname{acl}}
\nc{\dcl}{\operatorname{dcl}}
\nc{\tp}{\operatorname{tp}}
\nc{\stp}{\operatorname{stp}}
\nc{\Z}{\mathbb{Z}}
\nc{\N}{\mathbb{N}}
\nc{\M}{\mathfrak{M}}
\nc{\F}{\mathbb{F}}
\nc{\FF}{\mathfrak{F}}
\nc{\GG}{\mathfrak{G}}
\nc{\HH}{\mathfrak{H}}
\nc{\AS}{\wp}
\nc{\kappacap}{\bigcap\!{}^\kappa\,}

\def\Ind#1#2{#1\setbox0=\hbox{$#1x$}\kern\wd0\hbox to 0pt{\hss$#1\mid$\hss}
\lower.9\ht0\hbox to 0pt{\hss$#1\smile$\hss}\kern\wd0}
\def\Notind#1#2{#1\setbox0=\hbox{$#1x$}\kern\wd0\hbox to 0pt{\mathchardef
\nn="3236\hss$#1\nn$\kern1.4\wd0\hss}\hbox to 0pt{\hss$#1\mid$\hss}\lower.9\ht0
\hbox to 0pt{\hss$#1\smile$\hss}\kern\wd0}

\usepackage{babel}

\makeatother

\usepackage{babel}

\usepackage{babel}

\makeatother

\usepackage{babel}

\makeatother

\usepackage{babel}

\begin{document}
\global\long\global\long\global\long\def\alg{\operatorname{alg}}
\global\long\global\long\global\long\def\sep{\operatorname{sep}}
\global\long\global\long\global\long\def\K{\mathbf{K}}

\global\long\global\long\global\long\def\ins{\operatorname{ins}}
\global\long\global\long\global\long\def\dcl{\operatorname{dcl}}

\title{Artin-Schreier extensions in NIP and simple fields}

\author{Itay Kaplan, Thomas Scanlon and Frank O. Wagner}
\begin{abstract}
We show that NIP fields have no Artin-Schreier extension, and
that simple fields have only a finite number of them. 
\end{abstract}
\maketitle

\section{Introduction}

In \cite{mac}, Macintyre showed that an infinite $\omega$-stable
commutative field is algebraically closed; this was subsequently generalized
by Cherlin and Shelah to the superstable case; they also showed that
commutativity need not be assumed but follows \cite{chsh}. It is
known \cite{wood} that separably closed infinite fields are stable;
the converse has been conjectured early on \cite{Ch97}, but little
progress has been made. In 1999 the second author published on his
web page a note proving that an infinite stable field of characteristic
$p$ at least has no Artin-Schreier extensions, and hence no finite
Galois extension of degree divisible by $p$. This was later generalized
to fields without the independence property (NIP fields) by
(mainly) the first author.

In the simple case, the situation is even less satisfactory. It is
known that an infinite perfect, bounded (i.e.\ with only finitely
many extensions of each degree) PAC (pseudo-algebraically closed:
every absolutely irreducible variety has a rational point) field is
supersimple of SU-rank one \cite{hru}. Conversely, Pillay and Poizat
have shown that supersimple fields are perfect and bounded; it is
conjectured that they are PAC, but the existence of rational points
has only been shown for curves of genus zero (and more generally Brauer-Severi
varieties) \cite{spw}, certain elliptic or hyperelliptic curves \cite{mpp},
and abelian varieties over pro-cyclic fields \cite{mpw,il}. Bounded
PAC fields are simple \cite{ch} and again the converse is conjectured,
with even less of an idea on how to prove this \cite[Conjecture 5.6.15]{wa00};
note though that simple and PAC together imply boundedness \cite{ch99}.
In 2006 the third author adapted Scanlon's argument to the simple
case and showed that simple fields have only finitely many Artin-Schreier
extensions.

In this paper we present the proofs for the simple and the NIP
case, and moreover give a criterion for a valued field to be NIP
due to the first author.

We would like to thank Martin Hils and Fran\c coise Delon for some very
helpful comments and discussion on valued fields.

\section{Preliminaries}

\begin{notation}
\begin{enumerate}
\item Throughout the paper, $K$ will be a field of characteristic $p>0$. All fields considered will be supposed to be contained in a big algebraically closed field $\K$. We shall identify algebraic varieties defined over subfields of $\K$ with the set of their $\K$-rational points.
\item We denote by $K^{\alg}$ and $K^{\sep}$ the algebraic and separable closures of $K$, respectively. 
\item When we write $\bar{a}\in K$ for a tuple $\bar{a}=\left(a_{0},\ldots,a_{n}\right)$
we mean that $a_{i}\in K$ for $i\leq n$. 
\end{enumerate}
\end{notation}

\begin{defn}
A field extension $L/K$ is called an \emph{Artin-Schreier} extension
if $L=K(\alpha)$ for some $\alpha\in L\setminus K$ such that $\alpha^{p}-\alpha\in K$.
\end{defn}
Note that if $\alpha$ is a root of the polynomial $x^{p}-x-a$ then
$\{\alpha,\alpha+1,\ldots,\alpha+p-1\}$ are all the roots of the
polynomial. Hence, if $\alpha\notin K$ then $L/K$ is Galois and
cyclic of degree $p$. The converse is also true: if $L/K$ is Galois
and cyclic of degree $p$ then it is an Artin-Schreier extension \cite[Theorem VI.6.4]{lang}.

Let $\AS:K\rightarrow K$ the additive homomorphism given by $\AS(x)=x^{p}-x$.
Then the Artin-Schreier extensions of $K$ are bounded by the number
of cosets in $K/\AS(K)$. Indeed, since $K$ contains $\ker\AS=\mathbb{F}_{p}$,
if $K(\alpha)$ and $K(\beta)$ are two Artin-Schreier extensions,
then $a=\AS(\alpha)$ and $b=\AS(\beta)$ are both in $K\setminus\AS(K)$,
and \[
a-b=\AS(\alpha-\beta)\in\AS(K)\]
 implies $\alpha-\beta\in K$ and hence $K(\alpha)=K(\beta)$.
\begin{rem}
\label{rmk:ASbijection}In fact, the Artin-Schreier extensions of
$K$ are in bijection with the non-trivial orbits under the action of $\F_{p}^{\times}$
on $K/\AS(K)$. \end{rem}
\begin{proof}
Let $G=\mbox{Gal}(K^{\sep}/K)$. From \cite[X.3]{serre} we know that
$K/\AS(K)$ is isomorphic to $\mbox{Hom}(G,\Z/p\Z)$ as additive groups. This isomorphism respects the action of $\F_{p}^{\times}$. Hence the non-trivial orbits of this action on $K/\AS(K)$ are in bijection with the non-trivial orbits of the action on $\mbox{Hom}(G,\Z/p\Z)$. \\
Since $\mbox{Aut}(\Z/p\Z)=\F_p^\times$, these orbits are in bijection with \[\{\ker(\varphi) \mid 0 \neq \varphi \in\mbox{Hom}(G,\Z/p\Z)\}.\] By Galois theory, these are in bijection with Artin-Schreier extensions.
\end{proof}
We now turn to vector groups.

\begin{defn}
A \emph{vector group} is an algebraic group isomorphic as an algebraic
group to a finite Cartesian power of the additive group of a field. \end{defn}
\begin{fact}
\label{fac:Humphreys}\cite[20.4, Corollary]{humph} A closed connected
subgroup of a vector group is a vector group. 
\end{fact}
Using infinite Galois cohomology (namely, that $H^{1}\left(\mbox{Gal}(K^{\sep}/K),(K^{\sep})^{\times}\right)=1$,
for more on that see \cite[X]{serre}), one can deduce the following
fact: 
\begin{cor}
\label{fac:BlackBox} Suppose $K$ is perfect, and $G$ a closed connected
1-dimensional algebraic subgroup of $\left(\K,+\right)^{n}$ defined
over $K$, for some $n<\omega$. Then $G$ is isomorphic over $K$
to $\left(\K,+\right)$.
\end{cor}
This fact can also proved by combining Th\'eor\`eme 6.6 and Corollaire
6.9 in \cite[IV.3.6]{demazure}. Note that since $K$ is perfect, ``defined over $K$'' in the sense of algebraic geometry is the same as ``definable over $K$ in the structure $\K$'' in the model-theoretic sense.

We shall be working with the following group:
\begin{defn}
\label{defn} Let $(a_{1},\ldots,a_{n})=\bar{a}\in\K$. Put \[
G_{\bar{a}}=\{(t,x_{1},\ldots,x_{n})\in\K^{n+1}\mid t=a_{i}\cdot\AS\left(x_{i}\right))\mbox{ for }1\leq i\leq n\}.\]

\end{defn}
This is an algebraic subgroup of $(\K,+)^{n+1}$.

Recall that for an algebraic group $G$ we denote by $G^{0}$ the
connected component (subgroup) of the unit element of $G$. Note that
if $G$ is definable over $\bar{a}\in\K$, then $G^{0}$
is definable over $\bar{a}$ as well. Hence if $K$ is perfect
and $G$ is defined over $K$, then $G^{0}$ is defined over $K$.

\begin{lem}
\label{connected} If $\bar{a}\in\K$ is algebraically independent,
then $G_{\bar{a}}$ is connected. \end{lem}
\begin{proof}
By induction on $n:=\mbox{length}(\bar{a})$. If $n=1$, then $G_{\bar{a}}=\left\{ (t,x)\mid t=a_{1}\cdot\AS\left(x\right)\right\} $
is the graph of a morphism, hence isomorphic to $\left(\K,+\right)$
and thus connected. Assume the claim for $n$, and for some algebraically
independent $\bar{a}\in K$ of length $n+1$ let $\bar{a}'=\bar{a}\upharpoonright n$.
Consider the projection $\pi:G_{\bar{a}}\rightarrow G_{\bar{a}'}$.
Since $\K$ is algebraically closed, $\pi$ is surjective. Let $H=G_{\bar{a}}^{0}$
be the identity connected component of $G_{\bar{a}}$. As $[G_{\bar{a}'}:\pi(H)]\leq[G_{\bar{a}}:H]<\infty$,
it follows that $\pi(H)=G_{\bar{a}'}$ by the induction hypothesis.
Assume that $H\neq G_{\bar{a}}$.
\begin{claim}
For every $(t,\bar{x})\in G_{\bar{a}'}$ there is exactly one $x_{n+1}$
such that $(t,\bar{x},x_{n+1})\in H$. \end{claim}
\begin{proof}
Suppose for some $(t,\bar{x})$ there were $x_{n+1}^{1}\neq x_{n+1}^{2}$
such that $(t,\bar{x},x_{n+1}^{i})\in H$ for $i=1,2$. Hence their
difference $(0,\bar{0},\alpha)\in H$. But $0\neq\alpha\in\F_{p}$
by definition of $G_{\bar{a}}$. Hence, $(0,\bar{0},1)\in H$, and
$(0,\bar{0},\beta)\in H$ for all $\beta\in\mathbb{F}_{p}$. We know
that for every $(t,\bar{x},x_{n+1})\in G_{\bar{a}}$ there is some
${x'}_{n+1}$ such that $(t,\bar{x},{x'}_{n+1})\in H$; as $x_{n+1}-{x'}_{n+1}\in\F_{p}$
we get $(t,\bar{x},x_{n+1})\in H$ and $G_{\bar{a}}=H$, a contradiction.
\end{proof}
So $H$ is a graph of a function $f:G_{\bar{a}'}\rightarrow\K$ defined
over $\bar{a}$. Now put $t=1$ and choose $x_{i}\in\K$ for $i\leq n$
such that $a_{i}\cdot\AS\left(x_{i}\right)=1$. Let $L=\F_{p}(x_{1},\ldots,x_{n})$
and note that $a_{i}\in L$ for $i\leq n$. Then \[
x_{n+1}:=f(1,\bar{x})\in\dcl(a_{n+1},x_{1},\ldots,x_{n})=L(a_{n+1})_{\ins},\]
 where $L(a_{n+1})_{\ins}$ is the inseparable closure $\bigcup_{n<\omega}L(a_{n+1})^{p^{-n}}$
of $L(a_{n+1})$. Since $x_{n+1}$ is separable over $L(a_{n+1})$,
it follows that $x_{n+1}\in L(a_{n+1})$. By assumption, $a_{n+1}$
is transcendental over $\bar{a}'$, whence over $L$, and so $x_{n+1}\notin L$.
Hence $x_{n+1}=h(a_{n+1})/g(a_{n+1})$ for some mutually prime polynomials
$g,h\in L[X]$. But then \[
a_{n+1}\cdot\big[h(a_{n+1})^{p}/g(a_{n+1})^{p}-h(a_{n+1})/g(a_{n+1})\big]=1\]
 implies the equality in $L\left[X\right]$: \[
X\cdot\big[h^{p}-hg^{p-1}\big]=g^{p}.\]
 This implies that $h$ divides $g^{p}$, whence $h\in L$ is constant.
Similarly, $g$ divides $X$, which easily yields a contradiction.\end{proof}
\begin{cor}
\label{cla:G_bIsomorph}If $K$ is perfect and $\bar{a}\in K$
then $G_{\bar{a}}^{0}$ is isomorphic over $K$ to $(\K,+)$. In particular,
for any field $L\supseteq K$ with $L\subseteq\K$ the group $G_{\bar{a}}^{0}(L)$ is
isomorphic to $(L,+)$. In the case where $\bar{a}$ is algebraically
independent, $G_{\bar{a}}^{0}=G_{\bar{a}}$, so the same is true for
$G_{\bar{a}}$.\end{cor}
\begin{proof}
Over $\K$ the projection to the first coordinate of $G_{\bar{a}}$
is onto and has finite fibers, so $\dim G_{\bar{a}}=1$ (as a variety).
But then $G_{\bar{a}}^{0}$ is isomorphic over $K$ to $(\K,+)$ by
Corollary \ref{fac:BlackBox};
this isomorphism sends $G_{\bar{a}}^{0}(L)$ onto $(L,+)$. Finally,
if $\bar{a}$ is algebraically independent, $G_{\bar{a}}=G_{\bar{a}}^{0}$
by Lemma \ref{connected}. 
\end{proof}

\section{Simple fields}

For background on simplicity theory, the interested reader may consult
\cite{wa00}. The only property we shall need is a type-definable
variant of Schlichting's Theorem.
\begin{fact}
\label{schl}\cite[Theorem 4.5.13]{wa00} Let $G$ and $\Gamma$ be
type-definable groups with a definable action of $\Gamma$ on $G$,
and let $\FF$ be a type-definable $\Gamma$-invariant family of subgroups
of $G$. Then there is a $\Gamma$-invariant type-definable subgroup
$N\leq G$ containing some bounded intersection of groups in $\FF$
such that $[N:N\cap F]$ is bounded for all $F\in\FF$.\end{fact}
\begin{thm}
Let $K$ be a type-definable field in a simple theory. Then $K$ has
only boundedly many Artin-Schreier extensions.
\end{thm}
This means that in any elementary extension $\M$, the number of Artin-Schreier
extensions of $K^{\M}$ remains bounded. In particular, by compactness,
if $K$ is definable, it has only finitely many Artin-Schreier extensions.
\begin{proof}
If $K$ is finite, then it has precisely one Artin-Schreier extensions.
So we may assume it is infinite, and that the model is sufficiently
saturated. Let $k=K^{p^{\infty}}=\bigcap K^{p^{n}}$, a perfect infinite
type-definable sub-field. We shall show that $\AS(K)$ has bounded index in $K$ (as an additive subgroup).

Let $\FF=\{a\,\AS(K)\mid a\in k^\times\}$; this is a type-definable $k^{\times}$-invariant
family of additive subgroups of $K$. By Fact \ref{schl} there exists
a type-definable additive $k^{\times}$-invariant subgroup $N\leq K$
containing a bounded intersection of groups in $\FF$, say $\bigcap_{a\in A}a\,\AS(K)$ for some bounded $A\subset k$, such that $[N:N\cap F]$
is bounded for all $F\in\FF$.

For any finite $\bar{a}\in A$ the group $G_{\bar{a}}^{0}(k)$
is isomorphic to $(k,+)$ by Corollary \ref{cla:G_bIsomorph}. Now the projection of $G_{\bar{a}}(k)$ to the first coordinate is equal to $\bigcap_{a\in\bar{a}}a\,\AS(k)$; since the fibers of the projection are finite and $k$ is infinite, $\bigcap_{a\in\bar{a}}a\,\AS(k)$ must be infinite, too. By compactness $\bigcap_{a\in A}a\,\AS(k)$ is infinite, so $N\cap k$ is infinite as well. But $N\cap k$
is $k^{\times}$-invariant, hence an ideal in $k$, and must equal
$k$. Since $[N:N\cap\AS(K)]$ is bounded, so is $[k:k\cap\AS(K)]$.

Now $a=a^{p}+\AS(-a)$ for any $a\in K$, whence $K=K^{p}+\AS(K)$.
Assume $K=K^{p^{n}}+\AS(K)$. Then $K^{p}=K^{p^{n+1}}+\AS(K^{p})$,
whence \[
K=K^{p}+\AS(K)=K^{p^{n+1}}+\AS(K^{p})+\AS(K)=K^{p^{n+1}}+\AS(K);\]
 by compactness $K=k+\AS(K)$. Thus $[K:\AS(K)]=[k:k\cap\AS(K)]$
is bounded. \end{proof}
\begin{rem}
The important category of objects in simple theories are the hyper-definable
ones: Quotients of a type-definable set by a type-definable equivalence
relation. However, a hyper-definable field is easily seen to be type-definable:
If $K$ is given by a partial type $\pi$ modulo a type-definable
equivalence relation $E$, then for $a,b\in K$ the inequivalence$\neg aEb$
is given by the partial type $\exists x\,[\pi(x)\land(a-b)xE1]$.
By compactness, $E$ is definable on $\pi$. 
\end{rem}

\section{NIP fields}
\begin{defn}
A theory $T$ has the \emph{independence property} if there is a formula
$\phi(\bar{x},\bar{y})$ and some model $\M$ containing tuples $(\bar{a}_{i}:i\in\omega)$
and $(\bar{b}_{I}:I\subset\omega)$ such that $\M\models\phi(\bar{a}_{i},\bar{b}_{I})$
if and only if $i\in I$.

A theory $T$ is \emph{NIP} if it does not have the independence property. Such a theory is also called a dependent theory.\end{defn}
\begin{rem}
\label{additive} Let $K$ be an infinite field,
and let $f:K\to K$ be an additive polynomial, i.e. $f(x+y)=f(x)+f(y)$
for all $x,y\in K$. Then $f$ is of the form $\sum a_{i}x^{p^{i}}$.
Furthermore, if $K$ is algebraically closed and $\mbox{\ensuremath{\ker}}(f)=\F_{p}$,
then $f=a\cdot(x^{p}-x)^{p^{n}}$ for some $n<\omega$ and $a\in K$.\end{rem}
\begin{proof}
The first part appears in \cite[Proposition 1.1.5]{goss}. Assume
now that $K$ is algebraically closed and $|\ker(f)|=p$. If $a_{0}\neq0$,
then $(f,f')=1$, hence $f$ has no multiple factors and $deg(f)=p$.
If $a_{0}=0$, then $f=(g(x))^{p}$ for some additive polynomial $g$
with $|\ker(g)|=p$. So by induction $f=(a_{0}x+a_{1}x^{p})^{p^{n}}$
for some $n<\omega$. If moreover $\ker(f)=\F_{p}$, then $a_{0}+a_{1}=0$
hence $f=a\cdot(x^{p}-x)^{p^{n}}$ for some $a\in K$. \end{proof}
\begin{thm}
\label{main} Let $K$ be an infinite NIP field. Then $K$ is
Artin-Schreier closed.\end{thm}
\begin{proof}
We may assume that $K$ is $\aleph_{0}$-saturated, and we put $k=K^{p^{\infty}}$,
a type-definable infinite perfect sub-field (all
contained in an algebraically closed $\K$).
By dependence the Baldwin-Saxl condition \cite{bnsxl} holds, which
means that there is some $n<\omega$ such that for every $(n+1)$-tuple
$\bar{a}$ there is a sub-$n$-tuple $\bar{a}'$ with $\bigcap_{a\in\bar{a}}a\,\AS(K)=\bigcap_{a\in\bar{a}'}a\,\AS(K)$.
We fix some algebraically independent $(n+1)$-tuple $\bar{a}\in k$.
Let $G_{\bar a}$ be the group defined in \ref{defn}.
By Corollary \ref{cla:G_bIsomorph} we have algebraic isomorphisms
$G_{\bar{a}}\rightarrow(\K,+)$ and $G_{\bar{a}'}\rightarrow(\K,+)$
over $k$. Hence we can find an algebraic map $\rho$ over $k$ which
makes the following diagram commute: \[
\xymatrix{G_{\bar{a}}\ar[r]^{\pi}\ar[d] & G_{\bar{a}'}\ar[d]\\
(\K,+)\ar[r]^{\rho} & (\K,+)}
\]
As all groups and maps are defined over $k\subseteq K$, we can restrict
to $K$. By Baldwin-Saxl $\pi\upharpoonright G_{\bar{a}}(K)$ is onto $G_{\bar{a}'}(K)$,
so $\rho\upharpoonright K$ must be onto as well. Moreover, \[
|\ker(\rho)|=|\ker(\pi)|=|(0,\bar{0})\times\F_{p}|=p\,;\]
since $\ker(\pi)$ is contained in $G_{\bar{a}}(K)$, this remains
true in the restrictions to $K$. Finally, $\rho$ is a group homomorphism,
i.e.\ additive, and a polynomial, as it is an algebraic morphism
of $(\K,+)$.

Suppose that $0\neq c\in\ker(\rho)\subseteq k$, and put $\rho'(x)=\rho(c\cdot x)$.
Then $\rho'$ is an additive polynomial whose kernel is $\F_{p}$.
By Remark \ref{additive} there are $a\in \K$ (in
fact, $a\in k$ because $c\in k$, and $\rho$ is over $k$) and
$n<\omega$ such that $\rho'(x)=a\cdot(x^{p}-x)^{p^{n}}$. As $\rho'\upharpoonright K$
is onto $K$, for any $y\in K$ there is some $x\in K$ with \[
a\cdot(x^{p}-x)^{p^{n}}=a\cdot y^{p^{n}},\]
 so $\AS(x)=x^{p}-x=y$ and we are done.

In fact, $n$ must be $0$, as the degree of $\pi$ (as algebraic
morphism) is $p$, and so is the degree of $\rho'$, since the vertical
arrows are algebraic isomorphisms.\end{proof}
\begin{cor}
\label{extension}If $K$ is an infinite NIP field of characteristic
$p>0$ and $L/K$ is a finite separable extension, then $p$ does
not divide $[L:K]$. \end{cor}
\begin{proof}
Assume not, and let $L'$ be the normal closure of $L/K$. Then $p\mid[L':K]$,
so we may assume that $L/K$ is Galois. Let $G\leq\mbox{Gal}(L/K)$
be a subgroup of order $p$, and let $K^{G}\subseteq L$ be its fixed
field. As $K^{G}$ is interpretable in $K$, it is also NIP.
But $L/K^{G}$ is an Artin-Schreier extension, contradicting Theorem
\ref{main}. \end{proof}
\begin{cor}
\label{infinite} Let $K$ be an infinite NIP field of characteristic
$p>0$. Then $K$ contains $\F_{p}^{\alg}$. \end{cor}
\begin{proof}
Let $k=K\cap\F_{p}^{\alg}$, the relative algebraic closure of $\F_{p}$
in $K$. As $K$ is Artin-Schreier closed, so is $k$. Hence $k$
is infinite, perfect, and pseudo-algebraically closed. But \cite[Theorem 6.4]{duret}
of Duret states that a field with a relatively algebraically closed
PAC subfield which is not separably closed has the independence property.
Hence $k$ is algebraically closed, i.e.\
$k=\F_{p}^{\alg}$.
\end{proof}
One might wonder what happens for a type-definable field in a NIP
theory. We were unable to generalize our theorem to this case. However,
one easily sees:
\begin{prop}
\label{KPlusConnected}Let $K$ be a type-definable field in a NIP
theory. Then $K$ has either no, or unboundedly many Artin-Schreier
extensions. \end{prop}
\begin{proof}
By \cite{Shel} (another presentation appears in \cite[Proposition 6.1]{hpp})
there is a minimal type-definable subgroup $K^{00}$ of $(K,+)$ of
bounded index. As for any $\lambda\in K^{\times}$, the multiplicative
translate $\lambda K^{00}$ is also a type-definable additive subgroup
of bounded index, $K^{00}$ is an ideal of bounded index and must
therefore be equal to $K$. On the other hand, the image of $\AS$
is a type-definable subgroup of $(K,+)$. Remark \ref{rmk:ASbijection}
tells us that it has bounded index if and only if there are boundedly
many Artin-Schreier extensions. But if it has bounded index, then
it contains $K^{00}=K$, and $K$ is Artin-Schreier closed.
\end{proof}
We can, however, prove that a type-definable field is Artin-Schreier
closed under a stronger hypothesis.
\begin{defn}
Call a type-definable group $G$ \emph{strongly connected} when $G^{00}=G$. 
\end{defn}
Note that if $\pi:G\rightarrow H$ is a definable surjective group
homomorphism and $G$ is strongly connected, then so is $H$, since
$|G:\pi^{-1}(H^{00})|$ is bounded by $|H:H^{00}|$ and $\pi$ is
onto.
\begin{thm}
Let $K$ be a type-definable field in a NIP theory such that
there is no infinite decreasing sequence of type-definable additive
subgroups, each of unbounded index in its predecessor. Then $K$ is
Artin-Schreier closed.\end{thm}
\begin{proof}
We work in a saturated model. Let $\bar{a}=(a_{i}:i<\omega)$ be a
sequence of algebraically independent elements from $k=\bigcap K^{p^{n}}$.
Let $H_{i}=a_{i}\cdot\AS\left(K\right)$, and recall that $\bigcap_{j<n}H_{i_{j}}=\pi_{1}(G_{(a_{i_{0}},\ldots,a_{i_{n-1}})}(K))$
for all $i_{0}<\ldots<i_{n-1}$, where $\pi_{1}$ is the projection
to the first coordinate. Since $G_{(a_{i_{0}},\ldots,a_{i_{n-1}})}(K)$
is isomorphic (over $k$) to $(K,+)$ and we mentioned in \ref{KPlusConnected}
that the latter is strongly connected, $\bigcap_{j<n}H_{i_{j}}$ is
strongly connected, too. By assumption, there is some $n$ such that
$\bigcap_{i<n}H_{i}=\bigcap_{i<n+1}H_{i}$. Now proceed as in the
proof of Theorem \ref{main}.\end{proof}
\begin{rem}
\label{dep2}Saharon Shelah has shown \cite{Shel2} that this condition
holds when $T$ is strongly$^{2}$ dependent.
\end{rem}

\section{Some results on NIP valued fields}

Here we find a characterization of ``nice'' NIP valued fields
of characteristic $p>0$. First we recall the definitions and notations:
\begin{defn}
A valued field is a pair $\left(K,v\right)$ where $K$ is a field
and $v:K\rightarrow\Gamma\cup\left\{ \infty\right\} $ for an ordered
group $\Gamma$ such that: 
\begin{enumerate}
\item $v\left(x\right)=\infty$ if and only if $x=0$, 
\item $\Gamma=v\left(K^{\times}\right)$, 
\item $v\left(x\cdot y\right)=v\left(x\right)+v\left(y\right)$, and 
\item $v\left(x+y\right)\geq\min\left\{ v\left(x\right),v\left(y\right)\right\} $. 
\end{enumerate}
\end{defn}
If $\left(K,v\right)$ is a valued field, then $\Gamma$ is the \emph{valuation
group}, $\mathcal{O}_{K}=\{x\in K\mid v(x)\geq0\}$ is the (local)
ring of \emph{integers}, $\mathfrak{m}_{K}=\{x\in K\mid v(x)>0\}$
is its \emph{maximal ideal}, and $k=\mathcal{O}_{K}/\mathfrak{m}_{K}$
is the \emph{residue field}. As a structure we think of it as a 3-sorted
structure $(K,\Gamma,k)$ equipped with the valuation map $v:K^{\times}\rightarrow\Gamma$,
and the quotient map $\pi:\mathcal{O}_{K}\rightarrow k$. Other interpretations
are known to be equivalent (i.e.\ bi-interpretable, and hence to
preserve properties such as dependence).

In \cite{Del} Delon gave the following characterization of Henselian
NIP valued fields of characteristic $0$.
\begin{fact}
\cite{Del} Let $\left(K,v\right)$ be a Henselian valued field of
characteristic $0$. Then $\left(K,v\right)$ is NIP if and
only if the residue field $k$ is NIP. 
\end{fact}
Historically, this theorem stated that the valuation group must also
be NIP, but by a result of Gurevich and Schmitt \cite{GuSch},
every ordered abelian group is NIP.\\
Delon assumed initially that the field is radically closed, but later she saw that this assumption is redundant. In \cite{belair}, this fact appears as well, without this requirement.

Here we discuss valued fields of characteristic $p$, i.e.\ with
$\mathrm{char}(K)=\mathrm{char}(k)=p$.
\begin{prop}
\label{ResFieldInfi}If $(K,v)$ is a NIP valued field of characteristic
$p>0$, then the residue field contains $\F_{p}^{\alg}$.\end{prop}
\begin{proof}
Suppose $\pi\left(a\right)\in k$ with $a\in\mathcal{O}_{K}$. Since
$K$ is Artin-Schreier closed, there is $b\in K$ with $b^{p}-b=a$.
If $v(b)<0$, then $v(b^{p})=p\, v(b)<v(b)$, whence $v(b^{p}-b)=v(b^{p})<0$,
contradicting $v(a)\geq0$. Hence $v(b)\geq0$ and $b\in\mathcal{O}_{K}$.
Thus $\pi\left(b\right)\in k$, and $\pi\left(b\right)^{p}-\pi\left(b\right)=\pi\left(a\right)$.
In other words, $k$ is also Artin-Schreier closed, and hence infinite;
since it is interpretable, it is NIP, and contains $\F_{p}^{\alg}$
by Corollary \ref{infinite}.\end{proof}
\begin{prop}
\label{GroupPdiv}If $(K,v)$ is a NIP valued field of characteristic
$p>0$, then the valuation group $\Gamma$ is $p$-divisible.\end{prop}
\begin{proof}
Let $0>\alpha\in\Gamma$. So $\alpha=v\left(a\right)$ for some $a\in K^{\times}$.
As $K$ is Artin-Schreier closed, there is some $b\in K^{\times}$
such that $b^{p}-b=a$. Clearly $v(b)\geq0$ is impossible. Hence
$v(b^{p})=p\ v(b)<v(b)$, and \[
\alpha=v(a)=v(b^{p}-b)=\min\{v(b^{p}),v(b)\}=v(b^{p})=p\, v(b).\]
 So $\alpha$ is $p$-divisible, as is $\Gamma$ (for $\alpha$ positive,
consider $-\alpha$).
\end{proof}
As a corollary we obtain a result of Cherlin \cite{cher}. 
\begin{cor}
$\F_{p}\left(\left(t\right)\right)$ has the independence property, and so does $\F_{p}^{\alg}\left(\left(t\right)\right)$.
\end{cor}
Propositions \ref{ResFieldInfi} and \ref{GroupPdiv} are also sufficient
for a valued field to be NIP, under certain conditions. In order
to explain these conditions, we give two definitions.
\begin{defn}
A valued field $\left(K,v\right)$ of characteristic $p>0$ is called
a Kaplansky field if it satisfies: 
\begin{enumerate}
\item The valuation group $\Gamma$ is $p$-divisible, 
\item The residue field $k$ is perfect, and does not admit a finite separable
extension whose degree is divisible by $p$. 
\end{enumerate}
\end{defn}
This definition is taken from the unpublished book on valuation theory
by Franz-Viktor Kuhlmann \cite[13.11]{kuhlmann}. It is first-order
expressible, as the second condition is equivalent to saying that
for every additive polynomial $f\in k\left[x\right]$, and every $a\in k$,
there is a solution to $f\left(x\right)=a$ in $k$ (for a proof,
see \cite[Theorem 5]{KuhlmannAddPol}).
\begin{defn}
A valued field $\left(K,v\right)$ is called algebraically maximal
if it does not admit any non-trivial algebraic immediate extension
(i.e.\ keeping both the residue field and the valuation group). 
\end{defn}
This is also first order axiomatizable \cite[Chapter 14, Section 2]{kuhlmann}.
It always implies Henselianity, and is equivalent to it in characteristic
$0$. In characteristic $p$, it is weaker than being Henselian and
defectless (\cite[9.39]{kuhlmann}).

We shall use the following result of B\'elair.
\begin{fact}
\cite[Corollaire 7.6]{belair}\label{Belair} A valued field $K$
of characteristic $p$ which is Kaplansky and algebraically maximal
is NIP if and only if $k$ is NIP. 
\end{fact}
Finally, we have:
\begin{thm}
\label{VFChar}Let $(K,v)$ be an algebraically maximal valued field
of characteristic $p$ whose residue field $k$ is perfect. Then $(K,v)$
is NIP if and only if $k$ is NIP and infinite and $\Gamma$
is $p$-divisible.\end{thm}
\begin{proof}
If $\left(K,v\right)$ is NIP then $k$ is infinite (it even
contains $\F_{p}^{\alg}$), and NIP, and $\Gamma$ is $p$-divisible,
by Propositions \ref{ResFieldInfi} and \ref{GroupPdiv}. On the other
hand, if $k$ is NIP and infinite, by Corollary \ref{extension}
we get that $\left(K,v\right)$ is Kaplansky and we can apply fact
\ref{Belair}. 
\end{proof}
It is interesting to note the connection to Kuhlmann's notion of a
tame valued field (see \cite[Chapter 13, Section 9]{kuhlmann}). A
valued field $\left(K,v\right)$ is called tame if and only if it
is algebraically maximal, $\Gamma$ is $p$-divisible and $k$ is
perfect. Note the difference between this and Kaplansky.

We get as an immediate corollary:
\begin{cor}
Let $(K,v)$ be an algebraically maximal NIP valued field. Then
$K$ is tame if and only if $K$ is Kaplansky, if and only if $k$
is perfect.
\end{cor}

\end{document}